%Musterdokument: doc.tex
\input amstex
\documentstyle{amsppt}
\magnification=1200
%\input avier %Anweisungen fuer dina4-Format
%Es folgt die "Praeambel", falls erforderlich
\overfullrule=0pt
\NoRunningHeads
%
%\input amacro1
%%%Makros-Dokument%%%allgemeine Makros
%%%%%%%%%%%%%%%%%%%%%%%%%%%%%%%%%%%%%%%%%%%%%
%Blackbold Buchstaben C,F,H,K,N,P,R,Q,Z: z.B. \C
%%%%%%%%%%%%%%%%%%%%%%%%%%%%%%%%%%%%%%%%%%%%%
\def\C{\Bbb C}

\def\R{\Bbb R}

\def\Q{\Bbb Q}
\def\Z{\Bbb Z}
%%%%%%%%%%%%%%%%%%%%%%%%%%%%%%%%%%%%%%%%%%%%%%

%
\topmatter
\title
One remark on construction of separated quotient-space
\endtitle
%
%Titel wird fortlaufende Kopfzeile auf den ungeraden Seiten, wenn %nicht
%anders abgegeben
%\rightheadtext{...}
\author
Yu.A. Neretin
\endauthor
%
%
%Autor wird Kopfzeile auf den geraden Seiten, wenn nicht anders angegeben
%\leftheadtext{...}
%\affil 
%Was bedeutet affill?
%\endaffil
%
%
%\address
%Adresse
%\endaddress
\address
Adress Moscow State Institute of electronic and mathematics
\endaddress
%
%
%\curraddr
%derzeitige Adresse
%\endcurraddr
\curraddr
 Max-Planck-Institute fur Mathematik, Bonn
\endcurraddr
%%
%\email
%E-Mail-Adresse
%\endemail
\email
 neretin\@ mpim-bonn.mpg.de ; neretin\@ matan.miemstu.msk.su
\endemail
\date
Datum 31 december 1995
\enddate
%
%
%\dedicatory
%Widmung
%\enddedicatory
%
%
%\thanks
%Dank
%\endthanks
%
%
%\translator
%Uebersetzer
%\endtranslator
%
%
%
%
%\subjclass
%Klassifizierung
%\endsubjclass
%
%
\abstract
 We discuss elementary constructions of boundaries of symmetric spaces.
\endabstract
%\input toc...
%Muster fuer Inhaltsverzeichnis: toc.tex
%\toc
%oder
%\toc\nofrills{Eigene Ueberschrift}
%\widestnumber\specialhead{...}
%\widestnumber\head{...}
%\widestnumber\subhead{...}
%\widestnumber\subsubhead{...}
%\head {..} ... \page{..}
%\endhead
%
%\endtoc
\endtopmatter
\document
%\hyphenation{...}

Let $M$ be a compact metric space. Let
$M =
\underset{\alpha \in A} \to \cup M_{\alpha}$
be a partition of $M$
($M_{\alpha} \cap M_{\beta} = \phi$
if
$\alpha \ne \beta$).
Then the quotient-space $A$ has  canonical structure of a topological space . Recall that the set 
$P \subset A$
is closed if and only if
$\underset{\alpha \in P}\to \cup M_{\alpha}$
is a closed subset in $M$. Let
$a_1, a_2, \ldots$
be a sequence in $A$. Let
$a \in A$.
Then
$a_j \rightarrow a$
if there exist points
$m_j \in M_{\alpha_j}, m \in M_{\alpha}$
such that
$m_j \rightarrow m$
in $M$.

The space $A$ is not need to be separated in Hausdorff sence. We are interested in the following question: how to construct separated analog of the quotientspace $A$ ?

\head{1. Preliminaries. Hausdorff convergence} \endhead

Let
$N \subset M$
be a closed subset. Denote by
$M_{\varepsilon}$
the set of all points 
$m \in M$
satisfying the condition: there exist
$n \in N$
such that
$\rho (m,n) < \varepsilon$.
Let [M] be the space of all closed subsets in $M$. {\it Hausdorff distance}
$d(N,N^{\prime})$
in [M] between $N$ and
$N^{\prime}$
is the infimum of
$\varepsilon > 0$
such that
$N \subset N^{\prime}_{\varepsilon}$
and
$N^{\prime}_{\varepsilon} \subset N$.

Recall that the metric space [M] is compact. Recall also two simple facts on Hausdorff convergence. Denote by
$\bar S$
the closure of the set $S$. Denote by
$B_{\varepsilon} (m)$
the ball
$\rho (m,n) < \varepsilon$.

\demo{Lemma 1}
{\it Let
$N_j \in [M]$.
Let
$K_{\sigma} \ (\sigma \in \Sigma)$
be all limit points of the sequence
$N_j$.
Then
\roster
\item"a)" $\overline {\cup_{\sigma \in \Sigma} K}_{\sigma}$
coincides with the set of all
$m \in M$
such that for all
$\varepsilon > 0$
the set
$N_j \cap B_{\varepsilon} (m)$
is nonempty for infinite number of $j$.

\item"b)" $\cap_{\sigma \in \Sigma} K_{\sigma}$
coincides with the set of all
$m \in M$
such that for all
$\varepsilon > 0$
the set
$N_j \cap B_{\varepsilon} (m)$
is nonempty for sufficiently large $j$.
\endroster }
\enddemo

\head{2. Construction of separated quotient-space} \endhead

Let a partition
$M = \underset{\alpha \in A} \to \cup M_{\alpha}$
satisfies the following condition
\roster
\item"$\ast )$" for each
$B \subset A$
the set
$\overline{\underset {\alpha \in B}\to \cup M_{\alpha}}$
is the union of elements of the partition.
\endroster

 Fix an open subset
$\Cal A \subset A$
such that quotient-topology on
$\Cal A$
is separated.
Denote by
$\tilde \Cal A \subset [M]$
the set of subsets
$\overline M_{\alpha}, \, \alpha \in \Cal A$.
Let our data satify the condition

\roster 
\item "$\ast \ast $ )" the map
$\alpha \leftrightarrow \overline{M_{\alpha}}$
is a homeomorphism of the spaces
$\Cal A$
and
$\tilde \Cal A$.
\endroster

\demo{Definition}
The {\it sepapated quotient-space} [[A]] is the closure of
$\tilde \Cal A$
in Hausdorff metrics.
\enddemo

\demo{Remark}
Of course the construction depends on the set
$\Cal A \subset A$.
\enddemo

\head{3. Description of the set [[A]]} \endhead

By lemma 1 and the condition $\ast )$ the elements
$N \in [[A]]$
are unions of elements
$M_{\alpha}$
of the partition. Hence we associate to each
$N \in [[A]]$
subset
$S_N \subset A$
of all
$\sigma \in A$
such that
$M_{\alpha} \subset N$.
Denote by [A] the set of all subsets 
$S_N$.
By construction we have canonical bijection [[A]]
$\leftrightarrow$
[A].

The following proposition is evident.

\demo{Lemma 2}
{\it Let
$S \subset A$.
Then the following conditions are equivalent
\roster
\item"a)" $S \in [A]$
\item"b)" There exist a sequence
$a \in A$
such that each limit point of
$a_j$
is an element of $S$ and each element
$s \in S$ is a limit of the sequence
$a_j$
in yhe quotient-topology on $A$.
\endroster }
\enddemo

Elements of $[A]$ we call  {\it admissible subsets}. 

\head{4. Example: complete collineations} \endhead

Let $M$ be the Grassmann manifold
$Gr_n$
of all $n$-dimensional subspaces in
$\C^n \oplus \C^n$.

Let
$\lambda \in \C^* = \C \setminus 0$.
Let
$V \in Gr_n$.
Define the subspace $\lambda V$ :

$$h \oplus p \in V \Leftrightarrow h \oplus \lambda p \in \lambda V$$
where
$h \in \C^n \oplus 0, \, p \in 0 \oplus \C^n$.
Consider the partition of
$Gr_n$
into
$\C^*$
-orbits. Let
$Op \subset Gr_n$
be the space of graphs of invertible operators. Of course the space
$Op$
coincide with the general linear group
$GL_n(\C)$.
The quotient space
$Op/\C^* = GL_n(\C)/\C^*$
is the group
$PGL_n(\C)$
of invertible operators defined up to scalar multiplier.

We want to apply our construction to the space
$M = Gr_n$
and
$\Cal A = PGL_n(\C)$.
We have to describe all admissible subsets in
$Gr_n/\C^*$.

\demo{Example}
Let
$n=2$.
Consider the sequence
$Q_n = \pmatrix 1 & 0 \\ 0 & n \endpmatrix \in PGL_2(\C)$.
Then the set of limits of
$Q_n$
in
$Gr_2/\C^*$
consists of points
$V_1, \ldots, V_5$
(= subspaces in
$\C^2 \oplus \C^2$)
enumerated below:

$$\align &V_1 : (x,y; 0, 0) \\
&V_2 : (x, y; 0, y) \\
&V_3 : (x,0; 0, y) \\
&V_4 : (x, 0; x, y) \\
& V_5 : (0, 0; x,y) \endalign$$
where $x,y \in V$ . The subspaces
$V_1, V_3, V_5$
are stable points of the group
$\C^*$.
The $\C^*$-orbits of
$V_2, V_4$
are 1-dimensional complex curves.
\enddemo

\demo{Definition}
Let
$V \in Gr_n$.
Then
\roster
\item"a)" {\it Kernel}
$\text{Ker} \, V = V \cap (\C^n \oplus 0)$
\item"b)" {\it Image}
$\text{Im} \, V$
is the projection of $V$ to
$0 \oplus \C^n$.
\item"c)" {\it Domain}
$\text{Dom} \, V$
is the projection of $V$ to
$\C^n \oplus 0$.
\item"d)" {\it Indefiniteness}
$\text{Indef} \, V = V \cap (0 \oplus \C^n)$.
\endroster
\enddemo

\demo{Remark}
Let $V \in Gr_n$ Then the subspace $V$ induces by the 
obvious way the invertible operator

$$
\text{Dom} V/\text{Ker} V \rightarrow \text{Im V}/ \text{Indef} V
$$

We denote this operator by $<V>$.
\enddemo

\demo{Definition}
{\it Hinge} in
$\C^n$
is a collection

$$ {\Cal P}= (Q_0, P_1, Q_1, P_2, Q_2, \ldots, P_{k}, Q_{k})$$
where
$Q_j, P_j$
are elements of
$Gr_n$
defined up to multiplier and

0.
$$\align & Q_j = \text{Ker} \, Q_j \oplus \text{Indef} \, Q_j \\
& P_j \ne \text{Ker} \, P_j \oplus \text{Indef} \, P_j  \endalign$$
1. For each
$j = 1,2, \ldots, k$

$$\align & \text{Ker} \, P_j = \text{Ker} \, Q_j = \text{Dom} \, P_{j+1} \\
& \text{Im} \, P_j = \text{Im} \, Q_j = \text{Indef} \, P_{j+1} \endalign$$

2.
$$\align
&Q_0 = \C^n \oplus 0 \, \, ; \, \, \text{Dom} \, P_1 = \C^n \\
&Q_{k} = 0 \oplus \C^n  \, \, ; \, \, \text{Im} \, P_{k}= \C^n .  \endalign$$
\enddemo

\demo{Remark} 
Let $P$ be the graph of an invertible operator
$\C^n \rightarrow \C^n$.
Then

$$(\C^n \oplus 0, P, 0 \oplus \C^n)$$
is a hinge.
\enddemo

\demo{Remark}
The elements
$Q_0, \ldots, Q_{k +1}$
of a hinge are completely defined by the elements
$P_1, \ldots, P_{k}$. The subspaces
$Q_j$
are fixed points of the group
$\C^*$.
The
$\C^*$
-orbits of
$P_j$
are 1-dimensional complex curves.
\enddemo

\proclaim{Theorem}
The space
$[PGL_n]$
of all admissible subsets in
$Gr_n/\C^*$
coincides with the space of all hinges.
\endproclaim

The space
$[PGL_n]$
coincide with the {\it complete collineation} space constructed by Semple (see \cite{2}). It is a smooth algebraic variety and the group
$PGL_n$
is an open dense subset in
$[PGL_n]$.
On equivalence of these two constructions see see \cite{8}. Complete collineations is a partial case of complete symmetric varieties, see De Concini, Procesi  \cite{3}.

\head {5. Example. Furstenberg-Satake compactification of riemannian symmetric space} \endhead

 We will only discuss the case
$PGL_n(\R)/SO(n)$.
Consider the space
$\R^n \oplus \R^n$
provided by a skew-symmetric bilinear form
$\pmatrix 0 & 1 \\ -1 & 0 \endpmatrix$.
Let
$\Cal L$
be the grassmannian of all Lagrangian subspaces in
$\R^n \oplus \R^n$.
Denote by
$\R^*$
the multiplicative group of real positive numbers. This group acts on
$\Cal L$
by multiplications of linear relations on scalars.

Denote by $R$ the open subset in
$\Cal L$
consisting of graphs
of operators
$S : \R^n \rightarrow \R^n$.
It is easy to see that

$$\{\text{matrix $S$ is symmetric}\} \Leftrightarrow \{
\text{the grapf of $S$ is an element of $\Cal L$ }\}$$
The group
$GL_n(\R)$
acts on $R$ by the formula
$g : S \mapsto g^t Sg$.
The stabilizer of the point
$S = E$
is the orthogonal group
$O(n)$.
Hence $GL_n(\R)$-orbit $X$ of $E$ is a homogeneous space
$GL_n(\R)/O(n)$. Points of $X$ correpond to positive definite matrices $S$.

Now we apply the construction of the sections 2-3 to the space
$\Cal L$
and to the open subset
$X= GL_n(\R)/O(n)$.
Then the completion consists of hinges

$$P = (Q_0, P_1, Q_1, \ldots, P_{k}, Q_{k})$$
such that
$P_j \in \Cal L, Q_j \in \Cal L$ and the operators $<P_j>$
(see section 4) are positive definite.

\head {6. Example. Boundary of Bruhat-Tits building}\endhead

 Let
$Q_p$
be a $p$-adic field. Let $M$ be the space of all
$\Z_p$
-submodules in
$\Q_p$.
Let
$B \subset M$
be the space of all lattices. The group
$\Q^*_p$
act on $M$ in a natural way. Then the corresponding separated
quotient-space consists of collections

$$(R_0, T_1, R_1, \ldots, T_{k}, R_{k})$$
where
$0 = R_0 \subset T_1 \subset R_1 \subset T_2 \ldots \subset R_{k} = \Q^n_p$
are elements of $M$ defined up to multiplier,
$R_j$
are subspaces and images of
$T_j$
in
$R_j/R_{j-1}$
are lattices.

 I thanks C. De Concini, S.L.Tregub and E.B.Vinberg for discussion of this subject.

%\input biblio...
%Musterdokument fuer Literaturangaben: biblio.tex


\Refs\nofrills{Bibliography}

\ref
\no 1
\by E. Study
\paper \"Uber die Geometrie der Kegelschnitte, insbesondere deren charakteristische Probleme
\paperinfo Math. Ann., 27
\yr 1886
\pages 51-58
\endref

\ref
\no 2
\by I.G. Semple
\paper The variety whose points represent complete collineations of $S_r$ on $S^{\prime}_r$
\paperinfo Rend. Math. 10
\pages 201-280
\yr 1951
\endref

\ref
\no 3
\by F\"urstenberg, H.
\paper A Poisson formula for semisimple Lie groups
\paperinfo Ann. of Math. (2) 77
\yr 1963
\pages 335-386
\endref

\ref
\no 4
\by Satake, I.
\paper On representations and compactifications of symmetric Riemannian spaces
\paperinfo Ann. of Math. (2) 71
\yr 1960
\pages 77-110
\endref

\ref
\no 5
\by Alguneid, A.R.
\paper Complete quadrics primals in four dimensional space
\paperinfo Proc. Math. Phys. Soc. Egypt, 4
\yr 1952
\pages 93-104
\endref

\ref
\no 6
\by C. De Concini, C. Prochesi
\paper Complete symmetric varieties
\paperinfo Lect. Notes Math., 996
\pages 1-44
\endref

\ref
\no 7
\by Oshima, T., Sekiguchi, I.
\paper Eigenspaces of invariant differential operators on an affine symmetric space
\paperinfo Invent. Math. 57
\yr 1980
\pages 1-81
\endref

\ref
\no 8
\by Neretin, Yu. A.
\paper On universal completions of complex classical groups
\paperinfo Funct.Anal.Appl.,26:1
\pages 
\endref

\ref
\no 9
\by Vinberg E.B.
\paper On reductive algebraic semigroups
\paperinfo Adv.Sov.Math., volume dedicated to E.B.Dynkin, to appear
\endref

\ref
\no 10
\by Neretin Yu.A.
\paper Hinges and Study-Semple-Furstenberg-Satake-De Concini-Procesi-
Oshima boundary
\paperinfo to appear
\endref


\endRefs
\enddocument